\documentclass[11pt]{amsart}
\usepackage{amssymb}
\theoremstyle{plain}
\newtheorem{Thm}{Theorem}

\begin{document}

\title[XCF with boundary]
{Local existence to the cross curvature flow on 3-manifolds with
boundary}

\author{Li Ma, Baiyu Liu}

\begin{abstract}
In this paper, we use the DeTurck trick to study the short-time
existence of solutions to the Dirichlet and Newmann boundary
problems of the cross curvature flow on 3-manifolds with boundary.


\end{abstract}

 \maketitle

\section{Introduction}

In this paper, we study the Dirichlet and Newmann boundary
problems of the cross curvature flow on compact three-manifolds
with boundary. Let $(M, g)$ be a compact three-dimensional
Riemannian manifold with boundary $\partial M$ and with positive
sectional curvature or negative sectional curvature. The cross
curvature tensor $c=(c_{ij})$ on the Riemannian manifold $(M^3,g)$
is defined by
$$
c_{ij}\doteqdot det E(E^{-1})_{ij}=\frac{1}{2}\mu^{ipq}\mu^{jrs}E_{pr}E_{qs}
=\frac{1}{8}\mu^{pqk}\mu^{rsl}R_{ilpq}R_{kjrs},
$$
where $E_{ij}\doteqdot R_{ij}-\frac{1}{2}Rg_{ij}$ is the Einstein
tensor of the metric $g=(g_{ij})$, $det E\doteqdot \frac{det
E_{ij}}{det g_{ij}}$, and $\mu^{ijk}$ are the components of the
volume form $d\mu$ with indices raised.

We say that a 1-parameter family of metrics $(g(t)$ on 3-manifold
$M^3$ with the negative sectional curvature is a solution of the
cross curvature flow (XCF) if it satisfies
\begin{equation}\label{flow:nsc}
\frac{\partial }{\partial t}g=2c.
\end{equation}
Likewise, for the family $g(t)$ with the positive sectional
curvature, we say that $(M^3, g(t))$ is a solution if
\begin{equation}\label{flow:psc}
\frac{\partial }{\partial t}g=-2c.
\end{equation}

We have two local existence results.

\begin{Thm}\label{thm:der}
Let $(M^3, g_0)$ be a compact 3-manifold with boundary.  If the
sectional curvature of $(M^3, g_0)$ is either negative everywhere
or positive everywhere, then there is a unique short time solution
$g(t)$, $t\in [0, \epsilon)$, where $\epsilon>0$, to the cross
curvature flow (XCF) with
$$
g(x, 0)=g_0 (x),\quad x\in M
$$
and
$$
g(x, t)=g_0(x),\quad t\in [0,\epsilon) \quad x\in
\partial M.
$$
\end{Thm}

\begin{Thm}\label{thm:ste}
Let $(M^3, g_0)$ be a compact 3-manifold with boundary. Given any
smooth function $\lambda(t)$ depending only on the time variable
$t\in [0,+\infty)$. If the sectional curvature of $(M^3, g_0)$ is
either negative everywhere or positive everywhere, then there is a
unique short time solution $g(t)$, $t\in [0, \epsilon)$, where
$\epsilon>0$, to the cross curvature flow (XCF) with
\begin{equation}\label{con:inm}
g(x, 0)=g_0 (x),\quad x\in M
\end{equation}
and
\begin{equation}\label{con:umb}
h_{\alpha\beta}(x, t)=\lambda(t)g_{\alpha\beta}(x, t),\quad
\alpha, \beta=1,2, \quad x\in \partial M.
\end{equation}
\end{Thm}

We shall use the DeTurck trick to prove the above results. Since
the arguments for both results are similar, we shall only provide
the detail for the proof of Theorem \ref{thm:ste}.

The short-time existence of solutions to the cross curvature flow
on closed 3-manifolds has been obtained by B.Chow and Hamilton
\cite{CH} by the method of Nash-Moser implicit theorem. A simpler
proof for this result has been obtained by Buckland \cite{Buc05}
(see also \cite{Cho07}) and Buckland's argument used DeTurck
method \cite{Det83, Det03}. Interesting examples for the cross
curvature flows have been studied by X.Cao, Y.Ni, and Laurent
Saloff-Coste in \cite{Cao} and by L.Ma and D.Chen \cite{MC}.In the
case of curvature flows on manifolds of lower dimensions with
boundary, since the literature is huge, we only cite \cite{Bre02,
Bre02MathZ, Cor09, She96} here and one may find more references
therein. In \cite{She96}, Y. Shen has considered the Neumann
boundary value problem for the Ricci flow:
\begin{equation}\label{eq:rfb}
\left\{ \begin{array}{ll}
\frac{\partial g}{\partial t}=-2Ric, & x\in M\\
g(x, 0)=g_0, & x\in M\\
h_{\alpha \beta}(x, t)=\lambda g_{\alpha \beta}(x, t), & x\in
\partial M,
\end{array}\right .
\end{equation}
where $h_{\alpha \beta}$ is the second fundamental form of
$\partial M$ in $M$ and $\lambda$ is a constant. The short time
existence of (\ref{eq:rfb}) for such problem has been obtained in
\cite{She96}.

\section{Short Time Existence}

In this section, we present the proof of Theorem \ref{thm:ste} by
using the DeTurck trick.

Since the case of positive sectional curvature is similar,
we only consider the case of negative sectional curvature, namely
\begin{equation}\label{eq:xcf}
\left\{ \begin{array}{ll}
\frac{\partial g}{\partial t}=2c, & x\in M\\
g(x, 0)=g_0, & x\in M\\
h_{\alpha \beta}(x, t)=\lambda(t) g_{\alpha \beta}(x, t), & x\in \partial M.
\end{array}\right .
\end{equation}
We adopt the convention that Latin indices range from 1 to 3, while Greek indices range from 1 to 2.\\

\textbf{Step 1.} Analyze the linearization of (\ref{eq:xcf}).

We note that the linearization of (XCF) has been computed by
Buckland \cite{Buc05}. For reader's convenience, we give full
details.

If $\frac{\partial g_{ij}}{\partial s}=v_{ij} $ is a variation of the metric $g_{ij}$, then
\begin{equation}\label{eq:vco}
\frac{\partial }{\partial s}R_{ijkl}=\frac{1}{2}(\frac{\partial^2 v_{jl}}{\partial x^i \partial x^k}
+\frac{\partial^2 v_{ik}}{\partial x^j \partial x^l}-\frac{\partial^2 v_{jk}}{\partial x^i \partial x^l}-\frac{\partial^2 v_{il}}{\partial x^j \partial x^k})+\dots,
\end{equation}
where the dots denote the lower derivatives terms.
We normalize such that $\mu_{123}=\mu^{123}=1$ and recall from p.491 of \cite{Cho07} the following:
\begin{equation}\label{eq:mu}
\mu^{pqk}\mu^{rsl}R_{kjrs}=2E^{ml}(\delta_j ^p \delta_m ^q-\delta_m ^p \delta_j ^q).
\end{equation}
Applying (\ref{eq:mu}) to (\ref{eq:vco}), we obtain
\begin{eqnarray*}
  \frac{\partial}{\partial s}c_{ij}
  & = & \frac{1}{8}\mu^{pqk}\mu^{rsl}(\frac{\partial}{\partial
    s}R_{ilpq})R_{kjrs}+ \frac{1}{8}\mu^{pqk}\mu^{rsl} R_{ilpq}(\frac{\partial}{\partial s}R_{kjrs})\\
  & = & -\frac{1}{8}(\frac{\partial^{2} v_{lp}}{\partial x^{i}\partial x^{q}}
  +\frac{\partial^{2} v_{iq}}{\partial x^{l}\partial x^{p}}-\frac{\partial^{2} v_{lq}}{\partial x^{i}\partial x^{p}}
  -\frac{\partial^{2} v_{ip}}{\partial x^{l}\partial x^{q}})E^{ml}(\delta^{p}_{j}\delta^{q}_{m}-\delta^{p}_{m}\delta^{q}_{j})\\
  &   & +\frac{1}{8}(\frac{\partial^{2} v_{jr}}{\partial x^{k}\partial x^{s}}
  +\frac{\partial^{2} v_{ks}}{\partial x^{j}\partial x^{r}}-\frac{\partial^{2} v_{js}}{\partial x^{k}\partial x^{r}}
  -\frac{\partial^{2} v_{kr}}{\partial x^{j}\partial x^{s}}) E^{mk}(\delta^{r}_{i}\delta^{s}_{m}-\delta^{r}_{m}\delta^{s}_{i})+\dots\\
  & = & \frac{1}{2}E^{ml}(\frac{\partial^{2} v_{ij}}{\partial x^{l}\partial x^{m}}
  +\frac{\partial^{2} v_{lm}}{\partial x^{i}\partial x^{j}}-\frac{\partial^{2} v_{lj}}{\partial x^{i}\partial x^{m}}
  -\frac{\partial^{2} v_{im}}{\partial x^{l}\partial x^{j}})+\dots.
\end{eqnarray*}
Hence, the linearization of the map $X$ which takes $g$ to $2c$ is a second-order
partial differential operator. Its symbol is
\begin{equation}\label{eq:sox}
[\sigma DX(g)(\zeta )v]_{ij}=E^{ml}(\zeta_i\zeta_j v_{lm}
+\zeta_l\zeta_m v_{ij}-\zeta_i\zeta_m v_{lj}-\zeta_l\zeta_j v_{im}).
\end{equation}
Since this is homogenous, we may assume $\zeta$ has length one and rotate the
coordinates so that $\zeta=1$ and $\zeta_2=\zeta_3=0$. It follows that
$$
[\sigma DX(g)(\zeta )v]_{ij}=\delta_{i1}\delta_{j1}E^{ml}v_{lm}+v_{ij}E^{11}
-\delta_{i1}E^{1l}v_{lj}-E^{m1}\delta_{j1}v_{im}.
$$
We then deduce that
\begin{equation}\label{eq:sbx}\sigma DX(g)(\zeta)\left( \begin{array}{l}
v_{11}\\
v_{12}\\
v_{13}\\
v_{22}\\
v_{33}\\
v_{23}
\end{array}\right)
= \left( \begin{array}{llllll}
0 & 0 & 0 & E^{22} & E^{33} & 2E^{23}\\
0 & 0 & 0 & -E^{12} & 0 & -E^{13}\\
0 & 0 & 0 & 0 & -E^{13} & -E^{12}\\
0 & 0 & 0 & E^{11} & 0 & 0\\
0 & 0 & 0 & 0 & E^{11} & 0\\
0 & 0 & 0 & 0 & 0 & E^{11} \\
\end{array}\right)
\left( \begin{array}{l}
v_{11}\\
v_{12}\\
v_{13}\\
v_{22}\\
v_{33}\\
v_{23}
\end{array}\right)
.
\end{equation}
Notice that $E^{11}>0$, when the sectional curvature is negative. Therefore,
the eigenvalues of matrix $\sigma DX(g)(\zeta)$ is nonnegative, which indicates that (XCF) is weakly parabolic.\\

\textbf{Step 2.}
Modify the initial boundary value problem (\ref{eq:xcf}).

Following the idea from Shen \cite{She96}, we deduce a modified problem of (XCF).

Let $\tilde{g}$ be any background metric on $M$ with connection $\tilde{\Gamma}$.
Define the vector field such that its components are given by
\begin{equation}\label{eq:vfw}
W^k=g^{pq}(\Gamma^k_{pq}-\tilde {\Gamma}^k_{pq}).
\end{equation}
Now we are going to solve the following modified system:
\begin{equation}\label{eq:meq}
\left\{ \begin{array}{ll}
\frac{\partial g}{\partial t}=2c+L_W g, & x\in M\\
g(x, 0)=g_0, & x\in M\\
h_{\alpha \beta}(x, t)=\lambda(t)g_{\alpha \beta}(x, t), & x\in \partial M\\
g_{3\alpha}(x, t)=0, & x\in \partial M\\
W^3(x, t)=0, & x\in \partial M.
\end{array}\right .
\end{equation}

For a given point $x\in \partial M$, we can choose a local coordinate chart around $x$ such that $\{\frac{\partial }{\partial x^1},
\frac{\partial }{\partial x^2}\}$ form basis for $T_x\partial M$ and $\frac{\partial }{\partial x^3}$ is transversal to $T_x\partial M$. Therefore, one can check that
$$
e_3=\frac{g^{3i}}{(g^{33})^{\frac{1}{2}}}\frac{\partial }{\partial x^i},
$$
is the unit normal to $T_x\partial M$ in a small neighborhood of $x$.
The second fundamental form of $(\partial M, g_{\partial M} )$ in $(M, g)$ is
$$
h_{\alpha \beta}=-\frac{1}{2}L_{e_n}g_{\alpha \beta}=-\frac{1}{2}\frac{g^{3i}}{(g^{33})^{\frac{1}{2}}}\tilde{\nabla}_i g_{\alpha \beta}.
$$
The given boundary condition $g_{n\alpha}=0$ implies $g^{n\alpha}=0$ and
$$
h_{\alpha \beta}=-\frac{1}{2}(g^{33})^{\frac{1}{2}}\tilde{\nabla}_3 g_{\alpha \beta}.
$$
Combining the above equation with $h_{\alpha \beta}=\lambda(t)g_{\alpha \beta}$, we get
\begin{equation}\label{eq:sff}
\tilde{\nabla}_3 g_{\alpha \beta}=-\frac{2\lambda(t)}{(g^{33})^{\frac{1}{2}}}g_{\alpha \beta}.
\end{equation}

By the definition of $W$, we know
\begin{eqnarray*}
W^3 & = & \frac{1}{2}g^{ij}g^{3l}
(\tilde{\nabla}_jg_{3i}+\tilde{\nabla}_ig_{3j}-\tilde{\nabla}_lg_{ij})\\
& = & g^{ij}g^{3l}(\tilde{\nabla}_jg_{3i}-\frac{1}{2}\tilde{\nabla}_lg_{ij}).
\end{eqnarray*}
The given boundary condition $g_{3\alpha}=0$ implies $g^{3\alpha}=0$. Thus we have
\begin{eqnarray*}
W^3 & = & (g^{33})^2\tilde{\nabla}_3 g_{33}-\frac{1}{2}g^{33}g^{ij}g_{ij}\\
& = & \frac{1}{2}(g^{33})^2\tilde{\nabla}_3 g_{33}-\frac{1}{2}g^{33}g^{\alpha \beta}\tilde{\nabla }_3 g_{\alpha \beta}.
\end{eqnarray*}
Hence the condition $W^3=0$ is equivalent to
$$
\tilde{\nabla }_3 g_{33}=\frac{g^{\alpha \beta}\tilde{\nabla }_3 g_{\alpha \beta}}{g^{33}}
=-2\lambda(t)\frac{g^{\alpha \beta}g_{\alpha \beta}}{(g^{33})^{\frac{3}{2}}}=-4\lambda(t)\frac{1}{(g^{33})^{\frac{3}{2}}},
$$
where we have used (\ref{eq:sff}).

We conclude that (\ref{eq:meq}) is equivalent to
\begin{equation}\label{eq:mbc}
\left\{ \begin{array}{ll}
\frac{\partial g}{\partial t}=2c+L_W g, & x\in M\\
g(x, 0)=g_0, & x\in M\\
\tilde{\nabla}_3 g_{\alpha\beta}(x, t)=-\frac{2\lambda(t)}{(g^{33})^{\frac{1}{2}}}g_{\alpha \beta}(x, t), & x\in \partial M\\
g_{3\alpha}(x, t)=0, & x\in \partial M\\
\tilde{\nabla }_3 g_{33}=-4\lambda(t)\frac{1}{(g^{33})^{\frac{3}{2}}}, & x\in \partial M.
\end{array}\right .
\end{equation}

Next, we show that equation (\ref{eq:mbc}) has a short time solution, by proving that it is a parabolic equation.
The symbol $\sigma DX(g)(\zeta )v$ has been obtained in (\ref{eq:sox}). We now compute the symbol of the operator
$$
Y(g)\doteqdot L_W g.
$$
Compute
\begin{eqnarray*}
Y(g)_{ij} & = & \tilde{\nabla}_j W_{i}+\tilde{\nabla}_i W_{j}\\
& = & g_{ik}g^{pq}\tilde{\nabla}_j \Gamma^{k}_{pq}+g_{jk}g^{pq}\tilde{\nabla}_i \Gamma^{k}_{pq}\\
& = & \frac{1}{2}g_{ik}g^{pq}g^{kl}\frac{\partial }{\partial x^j}(\frac{\partial g_{pl}}{\partial x^q}
+\frac{\partial g_{ql}}{\partial x^p}-\frac{\partial g_{pq}}{\partial x^l})\\
&   & +\frac{1}{2}g_{jk}g^{pq}g^{kl}\frac{\partial }{\partial x^i}(\frac{\partial g_{pl}}{\partial x^q}+\frac{\partial g_{ql}}{\partial x^p}-\frac{\partial g_{pq}}{\partial x^l})+\dots\\
& = & \frac{1}{2}g^{pq}(\frac{\partial^2 g_{pi}}{\partial x^j\partial x^q}+\frac{\partial^2 g_{qi}}{\partial x^j\partial x^p}-\frac{\partial^2 g_{pq}}{\partial x^j\partial x^i})\\
&   & +\frac{1}{2}g^{pq}(\frac{\partial^2 g_{pj}}{\partial x^i\partial x^q}+\frac{\partial^2 g_{qj}}{\partial x^i\partial x^p}-\frac{\partial ^2 g_{pq}}{\partial x^i\partial x^j})+\dots.
\end{eqnarray*}
It follows that
\begin{eqnarray*}
\frac{\partial }{\partial s}Y(g)_{ij} & = & \frac{1}{2}g^{pq}(\frac{\partial^2 v_{pi}}{\partial x^j\partial x^q}
+\frac{\partial^2 v_{qi}}{\partial x^j\partial x^p}-\frac{\partial^2 v_{pq}}{\partial x^j\partial x^i})\\
&   & +\frac{1}{2}g^{pq}(\frac{\partial^2 v_{pj}}{\partial x^i\partial x^q}+\frac{\partial^2 v_{qj}}{\partial x^i\partial x^p}-\frac{\partial ^2 v_{pq}}{\partial x^i\partial x^j})+\dots\\
& = & g^{pq}(\frac{\partial^2 v_{pi}}{\partial x^j\partial x^q}+\frac{\partial^2 v_{pj}}{\partial x^i\partial x^q}-\frac{\partial^2 v_{pq}}{\partial x^j\partial x^i})+\dots,
\end{eqnarray*}
and
$$
[\sigma DY(g)(\zeta)v]_{ij}=g^{pq}(\zeta_j \zeta_q v_{pi}+\zeta_i \zeta_q v_{pj}-\zeta_j \zeta_i v_{pq}).
$$
Choosing an orthonormal basis and taking $\zeta_1=1$, $\zeta_2=\zeta_3=0$, we find that
\begin{equation}\label{eq:sol}\sigma DY(g)(\zeta)\left( \begin{array}{l}
v_{11}\\
v_{12}\\
v_{13}\\
v_{22}\\
v_{33}\\
v_{23}
\end{array}\right)
= \left( \begin{array}{llllll}
1 & 0 & 0 & -1 & -1 & 0\\
0 & 1 & 0 & 0 & 0 & 0\\
0 & 0 & 1 & 0 & 0 & 0\\
0 & 0 & 0 & 0 & 0 & 0\\
0 & 0 & 0 & 0 & 0 & 0\\
0 & 0 & 0 & 0 & 0 & 0 \\
\end{array}\right)
\left( \begin{array}{l}
v_{11}\\
v_{12}\\
v_{13}\\
v_{22}\\
v_{33}\\
v_{23}
\end{array}\right)
.
\end{equation}
Adding (\ref{eq:sbx}) and (\ref{eq:sol}) together, we finally arrive at
\begin{equation*}\sigma D(X+Y)(g)(\zeta)\left( \begin{array}{l}
v_{11}\\
v_{12}\\
v_{13}\\
v_{22}\\
v_{33}\\
v_{23}
\end{array}\right)
= \left( \begin{array}{llllll}
1 & 0 & 0 & E^{22}-1 & E^{33}-1 & 2E^{23}\\
0 & 1 & 0 & -E^{12} & 0 & -E^{13}\\
0 & 0 & 1 & 0 & -E^{13} & -E^{12}\\
0 & 0 & 0 & E^{11} & 0 & 0\\
0 & 0 & 0 & 0 & E^{11} & 0\\
0 & 0 & 0 & 0 & 0 & E^{11} \\
\end{array}\right)
\left( \begin{array}{l}
v_{11}\\
v_{12}\\
v_{13}\\
v_{22}\\
v_{33}\\
v_{23}
\end{array}\right)
,
\end{equation*}
which implies that (\ref{eq:mbc}) is parabolic.  Therefore, (\ref{eq:mbc})
has a unique short time solution, and so does (\ref{eq:meq}).
Let $g(x, t)$ $(t\in [0, \epsilon))$ be a solution of (\ref{eq:meq}).\\

\textbf{Step3.} Finally, we obtain the solution of (\ref{eq:xcf}) from the
solution of (\ref{eq:meq}), i.e. $g(x, t)$ $(t\in [0, \epsilon))$.\\
Let $\varphi_{t}(x)$ be the one-parameter family of diffeomorphism
which is determined by:
\begin{equation*}
\left\{ \begin{array}{ll}
\frac{\partial }{\partial t}\varphi_{t}(x)=W(\varphi(x), t), & x\in M\\
\varphi_{0}=id, & x\in M.
\end{array}\right .
\end{equation*}
Here $id$ denotes the identity map on M and vector field $W$ has components $W^k$ as defined by (\ref{eq:vfw}).
Let
$$
\bar{g}(y,t)=(\varphi_t ^{-1})^*(g(x, t)).
$$
We claim that
$\bar{g}(y,t)$ is the solution of (XCF) with conditions (\ref{con:inm}) and (\ref{con:umb}).
First, we compute
\begin{eqnarray*}
\frac{\partial }{\partial t}\bar{g}(y,t) & = & \frac{\partial }{\partial t}((\varphi_t ^{-1})^*(g(x, t)))\\
& = & (\varphi_t ^{-1})^*(\frac{\partial }{\partial t}g(x, t))-(\varphi_t ^{-1})^*(L_W g)\\
& = & (\varphi_t ^{-1})^*(2c(g)+L_W g)-(\varphi_t ^{-1})^*(L_W g)\\
& = & (\varphi_t ^{-1})^*(2c(g))\\
& = & 2c((\varphi_t ^{-1})^* g)\\
& = & 2c(\bar{g}).
\end{eqnarray*}
Secondly, we check that the boundary condition (\ref{con:umb}) is satisfied.
Let $\bar{h}$ be the second fundamental form of $(\partial M, \bar{g}|_{\partial M})$
in $(M, \bar{g})$, then we have
$$
\bar{h}_{\alpha \beta}=((\varphi_t ^{-1})^*(h))_{\alpha \beta}=\lambda(t)\bar{g}_{\alpha \beta}.
$$

Therefore, we have proved that (XCF) with condition
(\ref{con:inm}) and (\ref{con:umb}) has a unique short time
solution.

\end{document}